\documentclass[letterpaper,10pt,conference]{ieeeconf}
\IEEEoverridecommandlockouts
\usepackage{amsmath,amssymb}
\usepackage{mathtools}
\usepackage{graphicx}
\usepackage{cite}
\let\labelindent\relax
\usepackage{enumitem}

\newtheorem{thm}{Theorem}[section]
\newtheorem{lem}[thm]{Lemma}
\newtheorem{prp}[thm]{Proposition}
\newtheorem{ass}[thm]{Assumption}
\makeatletter
\@ifundefined{proof}{%
  \newenvironment{proof}{\par\noindent\textit{Proof:}\ }%
    {\hfill$\square$\par}%
}{}
\makeatother

\DeclareMathOperator{\tr}{tr}
\DeclareMathOperator{\logdet}{log\,det}
\DeclareMathOperator{\im}{im}
\DeclareMathOperator{\diag}{diag}
\newcommand{\norm}[1]{\lVert #1\rVert}
\newcommand{\st}{\text{s.t.}}

\newcommand{\mat}[1]{\begin{bmatrix}#1\end{bmatrix}}
\newcommand{\sqmat}[4]{\mat{#1&#2\\#3&#4}}
\newcommand{\matvec}[2]{\mat{#1\\#2}}
\newcommand{\hDel}{\hat{\Delta}}
\newcommand{\el}{\mathcal{E}}
\newcommand{\betat}[1]{\beta^{(#1)}}
\mathtoolsset{showonlyrefs}

\title{\LARGE \bf
Convex Reformulation of Minkowski-Sum Outer Approximation\\%
for Matrix Ellipsoids in Data-Driven Control
}
\author{Taira Kaminaga and Hampei Sasahara,
\IEEEmembership{Member, IEEE}
\thanks{This work was supported by JSPS KAKENHI Grant Number 24K17296,
by JST PRESTO Grant Number JPMJPR26K3
and JST SPRING Grant Number JPMJSP2108.
T. Kaminaga and H. Sasahara are with the Department of Information Physics and Computing, Graduate School of Information Science and Technology, the University of Tokyo, Tokyo, Japan \tt{\{tairakaminaga,hsasahara\}@g.ecc.u-tokyo.ac.jp}.}
}

\begin{document}
\maketitle
\thispagestyle{empty}
\pagestyle{empty}

\begin{abstract}

Approximating Minkowski sums of matrix ellipsoids by a single matrix ellipsoid provides a tractable description of noise bounds for controller design from noisy data. Existing linear matrix inequality (LMI) conditions can be used to compute such approximations, but directly solving the resulting optimization problems can be computationally expensive. This paper develops an efficient method for finding a minimum-volume outer matrix ellipsoid under these LMI conditions. Although the corresponding problem for ordinary vector-valued ellipsoids is convex, the objective function in the matrix-valued case is nonconvex. We find that the problem admits an equivalent convex reformulation, which resolves this difficulty. We first eliminate the matrix variables using equivalent LMI transformations and objective monotonicity, leaving one positive scalar variable per matrix ellipsoid in the sum. We then show that the objective becomes convex after a logarithmic change of variables. The resulting problem contains neither LMI constraints nor matrix variables and can therefore be solved efficiently. Numerical experiments confirm that the proposed method is substantially faster than the LMI-based method. An application to data-driven control further illustrates the tradeoff between guaranteed control performance and computation time.
\end{abstract}

\section{INTRODUCTION}

Ellipsoids are widely used to bound states and system parameters under bounded uncertainty~\cite{ellip:FOGEL_Automatica_1982,ellip:kurzhanski2000,ellip:durieu_multi-input_2001}. They have a compact representation in terms of a center and a shape matrix, and their linear images are easy to compute. These properties make them useful in state estimation and reachability analysis.
In these applications, state bounds are propagated through the dynamics and combined with disturbance bounds at each step. For linear dynamics and ellipsoidal bounds, this operation leads to Minkowski sums of ellipsoids~\cite{ellip:durieu_multi-input_2001,ellip:kurzhanskiy_TAC_2007}.
Since a Minkowski sum of ellipsoids is generally not an ellipsoid, an outer approximation is needed to retain an ellipsoidal representation. Computing a tight outer ellipsoid efficiently while guaranteeing containment is therefore a basic problem in these applications~\cite{ellip:durieu_multi-input_2001,ellip:kurzhanskiy_TAC_2007,ellip:yan_closed-form_2015,ellip:Halder_CDC_2018}.

Matrix ellipsoids extend this representation to matrix-valued uncertainty through quadratic matrix inequalities (QMIs). 
They are useful in data-driven control, where controllers for unknown systems are designed from noisy data. 
The measured data and a QMI encoding prior knowledge about the noise define a set of possible systems.
The goal is to find a controller that meets the desired specifications for every system in this set. 
When the system set and the control specifications admit QMI descriptions, the matrix S-procedure~\cite{DDCQMI:Waarde2022_TAC_origin,DDCQMI:Waarde2023_siam_qmi} yields linear matrix inequality (LMI) conditions for controller design. 
This approach has led to stabilization and performance guarantees in the presence of exogenous disturbances or measurement noise~\cite{DDCQMI:Waarde2022_TAC_origin,DDCQMI:Waarde2023_siam_qmi,DDCQMI:Steentjes2022_Cont_Sys_Let_Covariance,DDCQMI:Waarde2024_TAC_AR,DDCQMI:Kaminaga2025_ACC,DDCQMI:kaminaga2025datainformativity}.
This robust design approach has also been shown to be effective against bounded malicious perturbations of the data used for controller design~\cite{security:burbano_sasahara_BADControl,security:SASAHARA2026}.

A single matrix ellipsoid, however, may not adequately capture prior knowledge about the noise, such as instantaneous or element-wise bounds. It may also be insufficient to describe the superposition of multiple noise sources.
Several studies have considered such noise models in controller design~\cite{DDCQMI:BISOFFI2021,DDCQMI:BISOFFI2024_CSL,DDCQMI:Hu2025_RDPC}, and
recent work has introduced the structured data perturbation framework to describe these models in a unified way~\cite{DDCQMI:kaminaga2025datainformativity}. 
This framework represents the admissible perturbations as a Minkowski sum of linearly transformed QMI sets and provides sufficient LMI conditions for enclosing the sum in a single QMI set.
These LMIs can be incorporated into controller design to optimize the approximation and controller jointly. 
This method can handle various noise models, but its computational cost increases with the data length and the number of sets in the sum: both the number of decision variables and the LMI size grow. 


This motivates the development of an efficient method for computing tight outer approximations of the Minkowski sum. Volume is a natural measure of the tightness of the outer approximation. 
We seek a minimum-volume outer matrix ellipsoid subject to the LMI conditions in~\cite{DDCQMI:kaminaga2025datainformativity}. 
For ordinary vector-valued ellipsoids, the corresponding problem admits a convex formulation based on the S-procedure~\cite{boyd1994linear}.
In the matrix-valued case, the volume depends on two shape matrices.
Although the sufficient containment constraints can be expressed as LMIs, the resulting objective is generally nonconvex in the variables of this formulation.

We uncover a hidden convex structure in this LMI-constrained volume minimization problem. This leads to an equivalent convex reformulation with one scalar variable for each matrix ellipsoid in the sum without LMI constraints.
We first derive an equivalent form of the LMI conditions and show that the center and shape matrices of the optimal approximating matrix ellipsoid can be described in terms of positive scalar parameters.
This eliminates the matrix variables and LMI constraints. 
A logarithmic change of the parameters then makes the objective convex. 
We can thus solve volume minimization under the original LMI conditions through convex optimization over scalar variables.

We also prove that the objective is strictly convex on the normalized subspace, and derive explicit formulas for its gradient and Hessian. 
These formulas allow us to apply various convex optimization algorithms, including gradient descent method and Newton's method.
To assess the computational benefit of the reformulation, we evaluate the computation cost of each iteration in Newton's method, and compare the cost with that of directly handling the LMI in~\cite{DDCQMI:kaminaga2025datainformativity}.

In numerical experiments, the proposed method computes outer matrix ellipsoids with smaller volumes in substantially less time than an LMI-based baseline that restricts the approximation shape. We also apply the proposed formulation to data-driven control in two stages: first compute the outer approximation, then design the controller. 
In the tested settings, this approach gives weaker performance guarantees than the LMI-based method in~\cite{DDCQMI:kaminaga2025datainformativity}, but substantially reduces computation time. 
It thus offers an option for controller design when computational cost is a priority.

\subsection*{Related Work}

Research on Minkowski sums of ordinary ellipsoids includes reachable-set approximation~\cite{ellip:kurzhanski2000,ellip:kurzhanskiy_TAC_2007}, analytical descriptions of the boundary~\cite{ellip:yan_closed-form_2015}, and parameterizations of outer ellipsoids for optimization.
The study in~\cite{ellip:durieu_multi-input_2001} uses support functions to construct a family of outer ellipsoids with one scalar parameter for each ellipsoid in the sum. 
It also proves that volume minimization over this family is convex. The study in~\cite{ellip:Halder_CDC_2018} clarifies the relationships among existing parameterizations. Our parameterization reduces to that of~\cite{ellip:durieu_multi-input_2001} in the vector-valued case and extends it to matrix ellipsoids.

Another approach uses the S-procedure to derive LMI conditions for an ellipsoid to contain a Minkowski sum~\cite{boyd1994linear}. 
Volume minimization under these LMI conditions is also convex for ordinary ellipsoids.
The work in~\cite{DDCQMI:kaminaga2025datainformativity} extends these conditions to matrix variables and provides the starting point for our analysis. 
We characterize the optimal solutions of volume minimization under these conditions and show that the problem is equivalent to optimization over a parameterized family. 
Our results extend the parameterization in~\cite{ellip:durieu_multi-input_2001} and the LMI approach in~\cite{boyd1994linear} to the matrix-valued setting and establish their equivalence for volume minimization.

\subsection*{Notation and Organization}
$\mathbb{R}^n$, $\mathbb{R}_{++}^n$, and $\mathbb{R}_{+}^n$ denote the sets of real, strictly positive, and nonnegative $n$-vectors, respectively.
$\mathbb{S}^n$ denotes the set of real symmetric $n\times n$ matrices, and $\mathbb{S}_{++}^n$ its positive definite subset.
For a matrix $M$, $\im M$ denotes the image of $M$.
For symmetric $M$, $M\succ0$ ($M\succeq0$) and $M\prec0$ ($M\preceq0$) indicate positive and negative (semi)definiteness, respectively.
$A^{1/2}$ denotes the positive semidefinite square root of $A\succeq0$.
$I_n$ and $0_{n,m}$ denote the $n\times n$ identity matrix and the $n\times m$ zero matrix, respectively, with subscripts omitted when the dimensions are clear from context.
$\mathbf{1}$ denotes the vector of ones.
$\diag(A_1,\dots,A_K)$ or $\diag_{k=1}^K(A_k)$ denotes the block diagonal matrix with blocks $A_1,\dots,A_K$.
$A\otimes B$ denotes the Kronecker product of $A$ and $B$, $\bigoplus_{k=1}^K\mathcal{S}_k$ denotes the Minkowski sum of $\mathcal{S}_1,\dots,\mathcal{S}_K$, and $\nabla f(x)$ denotes the gradient of $f$ at $x$.

Section~\ref{sec:pre} reviews matrix ellipsoids and their role in data-driven control. Section~\ref{sec:prob} formulates the outer approximation problem and introduces a relaxation based on the matrix S-procedure. Section~\ref{sec:parametrize} reduces volume minimization to an optimization problem in scalar parameters. Section~\ref{sec:opt} derives the convex formulation and analyzes its properties and computational complexity. Section~\ref{sec:numerical} presents numerical experiments, and Section~\ref{sec:conclusion} concludes the paper.

\section{PRELIMINARIES}\label{sec:pre}
\subsection{Matrix Ellipsoids and Their Volume}\label{subsec:matellip}

We introduce matrix ellipsoids and summarize their basic properties.
Given $c\in\mathbb{R}^q$, $Q\in\mathbb{S}_{++}^q$, and $\sigma>0$, an ordinary ellipsoid is defined as
\begin{equation}\label{e:vecEllip}
    \{x\in\mathbb{R}^q\mid (x-c)^\top Q^{-1}(x-c)\le\sigma\},
\end{equation}
where $c$ is the center and $(Q,\sigma)$ determines the shape.
Since $(c,Q,\sigma)$ and $(c,\sigma Q,1)$ define the same ellipsoid, $\sigma$ can be normalized to $1$ without loss of generality.
We extend this definition to a matrix variable by replacing the scalar quadratic inequality with a QMI. The resulting matrix ellipsoid is
\begin{equation}\label{e:EllipDfn}
    \el(C,Q,R)\coloneqq\{X\in\mathbb{R}^{q\times r}\mid{}(X-C)^\top Q^{-1}(X-C)\preceq R\}.
\end{equation}
Here, $C\in\mathbb{R}^{q\times r}$ is the center, and the pair $(Q,R)\in\mathbb{S}_{++}^q\times\mathbb{S}_{++}^r$ determines the shape.
The ordinary ellipsoid in~\eqref{e:vecEllip} is recovered by setting $r=1$, $C=c$, and $R=\sigma$.
As in the vector-valued case, the representation in~\eqref{e:EllipDfn} has a scaling freedom in $Q$ and $R$:
\begin{equation}\label{e:scale_invariant}
    \el(C,Q/a,aR)=\el(C,Q,R),\qquad a>0.
\end{equation}
This freedom permits only scalar rescaling of $R$.
Fixing $R$ to a prescribed matrix, such as $I_r$, generally restricts the class of matrix ellipsoids that can be represented.

We define volume by identifying $\mathbb{R}^{q\times r}$ with $\mathbb{R}^{qr}$. The volume of $\el(C,Q,R)$ is then proportional to $\sqrt{\det(R\otimes Q)}$, with a proportionality constant that depends only on $q$ and $r$~\cite{DDCQMI:BISOFFI2021}. 
We therefore use $\Psi(Q,R)\coloneqq\logdet(R\otimes Q)$ as a measure of size.
Minimizing this quantity is equivalent to minimizing volume.

In data-driven control~\cite{DDCQMI:Waarde2022_TAC_origin,DDCQMI:Waarde2023_siam_qmi,DDCQMI:Steentjes2022_Cont_Sys_Let_Covariance,DDCQMI:Waarde2024_TAC_AR,DDCQMI:Kaminaga2025_ACC,DDCQMI:kaminaga2025datainformativity}, QMIs are written as
\begin{equation}\label{e:QMI}
    \matvec{I_r}{X}^{\top}\Phi\matvec{I_r}{X}\succeq0,
\end{equation}
where $\Phi\in\mathbb{S}^{q+r}$ is a coefficient matrix. This form connects the matrix ellipsoid representation to the formulation in~\cite{DDCQMI:kaminaga2025datainformativity}. To express~\eqref{e:EllipDfn} in the form~\eqref{e:QMI}, set
\begin{equation}
    \Phi=\Pi(C,Q,R)\coloneqq
    \sqmat{R-C^\top Q^{-1}C}{C^\top Q^{-1}}{Q^{-1}C}{-Q^{-1}}.
\end{equation}
Then $\el(C,Q,R)=\{X\mid\eqref{e:QMI},\ \Phi=\Pi(C,Q,R)\}$.

\subsection{Data-Driven Control Using Matrix Ellipsoids}\label{subsec:cnt}

The following control problem motivates our study of matrix ellipsoids. Consider a discrete-time linear time-invariant system $x_+=A^*x+B^*u$, where $x_+\in\mathbb{R}^n$, $x\in\mathbb{R}^n$, and $u\in\mathbb{R}^m$ denote the next state, current state, and input, respectively. The true system matrices $(A^*,B^*)$ are unknown, but noisy input-state data are available. In the absence of noise, offline data $(X_+^*,X^*,U^*)$ of length $T$ satisfy
\begin{equation}\label{cnt:cleandata}
    X_+^*=A^*X^*+B^*U^*,
\end{equation}
where $X_+^*\in\mathbb{R}^{n\times T}$, $X^*\in\mathbb{R}^{n\times T}$, and $U^*\in\mathbb{R}^{m\times T}$. To account for noise and disturbances in the available measurements $(X_+,X,U)$, we introduce additive perturbations $(\Delta_{X_+},\Delta_X,\Delta_U)$ by
\begin{equation}\label{cnt:measured_data}
    X_+=X_+^*+\Delta_{X_+},\quad X=X^*+\Delta_X,\quad U=U^*+\Delta_U.
\end{equation}
The perturbations have the same dimensions as the corresponding noise-free data matrices. Define the stacked perturbation matrix $\Delta=[\Delta_{X_+}^\top\ {-\Delta_X^\top}\ {-\Delta_U^\top}]^\top$, and assume that
\begin{equation}\label{cnt:noise_matellip}
    \Delta=E\hDel,\qquad\hDel^\top\in\el(C,Q,R),
\end{equation}
where $\hDel\in\mathbb{R}^{\hat n\times T}$. The matrices $E\in\mathbb{R}^{(2n+m)\times\hat n}$ and $(C,Q,R)\in\mathbb{R}^{T\times\hat n}\times\mathbb{S}_{++}^{T}\times\mathbb{S}_{++}^{\hat n}$ are assumed to be known. 
We seek a state-feedback controller $u=Kx$ that guarantees the desired closed-loop performance for every system consistent with the data $(X_+,X,U)$ and the prior knowledge of the perturbations. This setting includes several widely studied noise models.
The choice $E=[I_n\ 0_{n,n+m}]^\top$ corresponds to exogenous disturbances affecting only the next-state data~\cite{DDCQMI:Waarde2022_TAC_origin}.
Choosing $E=I_{2n+m}$ allows measurement noise in both the input and state data~\cite{DDCQMI:BISOFFI2024_CSL}.

The constraint~\eqref{cnt:noise_matellip} can represent an energy bound on the noise~\cite{DDCQMI:Waarde2022_TAC_origin}. A single matrix ellipsoid, however, often cannot capture more specific prior information. The structured data perturbation framework in~\cite{DDCQMI:kaminaga2025datainformativity} addresses this limitation by generalizing~\eqref{cnt:noise_matellip} to
\begin{equation}\label{cnt:noise_minkowski}
    \hDel^\top\in\bigoplus_{k=1}^{K}F_k\el(C_k,Q_k,R_k)G_k,
\end{equation}
where the known matrices $F_k$, $G_k$, $C_k$, $Q_k$, and $R_k$ encode prior information about the noise. This model includes the following cases.
\begin{enumerate}[label=(\roman*)]
\item \textit{Superposition of exogenous disturbance and measurement noise:}\label{cnt:superposition} 
Exogenous disturbances and measurement noise may be subject to separate bounds, each represented by a matrix ellipsoid. The set of combined perturbations is then the Minkowski sum of the linear images of these ellipsoids.
\item \textit{Instantaneously bounded perturbation:}\label{cnt:instantaneous_bound}
The noise vector at each time step satisfies a prescribed norm bound.
Such bounds are considered in~\cite{DDCQMI:BISOFFI2021,DDCQMI:BISOFFI2024_CSL}.
\item \textit{Element-wise bounded perturbation:}\label{cnt:elementwise_bound}
Each scalar noise component is constrained by known lower and upper bounds at each time step.
These bounds define an interval, which is a one-dimensional ellipsoid.
This model includes a common amplitude bound for all components and time steps, expressed as $\max_{i,j}|\hDel_{i,j}|\le\epsilon$.
\item These constraints can also be combined. For example,~\cite{DDCQMI:Hu2025_RDPC} considers exogenous disturbances and measurement noise with separate bounds on their absolute values. This setting combines~\ref{cnt:superposition} and~\ref{cnt:elementwise_bound}.
\end{enumerate}

To use such prior information in controller design, the method in~\cite{DDCQMI:kaminaga2025datainformativity} encloses the Minkowski sum in a single matrix ellipsoid. 
This allows control methods such as those in~\cite{DDCQMI:Waarde2022_TAC_origin} to be applied. The outer approximation enlarges the admissible noise set and, consequently, the set of systems consistent with the data. 
A controller that meets the specifications for every system in the enlarged set therefore also meets them for the true system. 
The method in~\cite{DDCQMI:kaminaga2025datainformativity} solves a sufficient LMI condition for containment together with the controller design conditions. 
In cases~\ref{cnt:instantaneous_bound} and~\ref{cnt:elementwise_bound}, however, its computational cost grows substantially with the data length. We therefore seek an efficient way to minimize the volume of the outer approximation independently of controller design.

\section{PROBLEM FORMULATION AND RELAXATION}\label{sec:prob}
\subsection{Problem Setting and Assumptions}

We consider the outer approximation of a Minkowski sum of matrix ellipsoids by a single matrix ellipsoid. Let $\el_k\coloneqq\el(C_k,Q_k,R_k)$, $k=1,\dots,K$, where $(C_k,Q_k,R_k)\in\mathbb{R}^{q_k\times r_k}\times\mathbb{S}_{++}^{q_k}\times\mathbb{S}_{++}^{r_k}$. Define the Minkowski sum of their linear images as
\begin{equation}
    \mathcal{T}\coloneqq\bigoplus_{k=1}^{K}F_k\el_kG_k,
\end{equation}
where $F_k\in\mathbb{R}^{q\times q_k}$ and $G_k\in\mathbb{R}^{r_k\times r}$ are known constant matrices for $k=1,\dots,K$. Our objective is to find a minimum-volume matrix ellipsoid that contains $\mathcal{T}$, that is, to solve
\begin{equation}\label{prb:prb}
    \min_{C\in\mathbb{R}^{q\times r},\ Q\in\mathbb{S}_{++}^q,\ R\in\mathbb{S}_{++}^r}
        \Psi(Q,R)
    \quad\st\quad \el(C,Q,R)\supseteq\mathcal{T}.
\end{equation}

Throughout the analysis, we make the following assumption.
\begin{ass}\label{ass:volume}
The matrices $F_k$ and $G_k$ are nonzero for every $k$. Furthermore, $\mathcal{T}$ has a nonempty interior.
\end{ass}

If either $F_k$ or $G_k$ is zero, then $F_k\el_kG_k=\{0\}$, so we can omit this term from $\mathcal{T}$. The nonempty interior assumption gives a positive lower bound on the volume of any outer matrix ellipsoid and implies the following rank conditions.
\begin{lem}\label{lem:FG}
Suppose that Assumption~\ref{ass:volume} holds. Then $F=[F_1\cdots F_K]$ has full row rank, and $G=[G_1^\top\cdots G_K^\top]^\top$ has full column rank.
\end{lem}
\begin{proof}
Suppose that $F$ does not have full row rank. Then there exists a nonzero vector $y\in\mathbb{R}^q$ such that $y^\top F=0$. Hence, $y^\top F_k=0$ for every $k$, so $y^\top X=0$ for every $X\in\mathcal{T}$. This contradicts the assumption that $\mathcal{T}$ has a nonempty interior. Thus, $F$ has full row rank. The argument for $G$ is analogous.
\end{proof}

\subsection{Approximation Based on the Matrix S-Procedure}

The containment condition $\el(C,Q,R)\supseteq\mathcal{T}$ is difficult to impose directly in numerical optimization. For ordinary ellipsoids, existing approaches use parameterized families of outer ellipsoids~\cite{ellip:durieu_multi-input_2001,ellip:kurzhanskiy_TAC_2007} or relaxations based on the S-procedure~\cite{boyd1994linear}. The work in~\cite{DDCQMI:kaminaga2025datainformativity} extends the latter approach to matrix ellipsoids. 
We use the resulting relaxation as the starting point of our analysis and state it below. 

\begin{prp}[{\cite[Theorem 6.3]{DDCQMI:kaminaga2025datainformativity}}]\label{prp:LMI}
Let $C\in\mathbb{R}^{q\times r}$, $Q\in\mathbb{S}_{++}^q$, and $R\in\mathbb{S}_{++}^r$, and set $\Phi=\Pi(C,Q,R)$ and $\Phi_k=\Pi(C_k,Q_k,R_k)$ for $k=1,\dots,K$. If there exists $\alpha\in\mathbb{R}_+^K$ such that
\begin{equation}\label{apx:LMI}
    \sqmat{I_r}{0}{0}{F^\top}\!\Phi\!\sqmat{I_r}{0}{0}{F}
    \succeq\sum_{k=1}^{K}\alpha_k\!
    \sqmat{G_k^\top}{0}{0}{U_k^\top}\!\Phi_k\!\sqmat{G_k}{0}{0}{U_k},
\end{equation}
then $\el(C,Q,R)\supseteq\mathcal{T}$. Here, $F=[F_1\cdots F_K]$ and
\begin{equation}
    U_k=\begin{bmatrix}
    0_{q_k,\sum_{l=1}^{k-1}q_l}&I_{q_k}&0_{q_k,\sum_{l=k+1}^{K}q_l}
    \end{bmatrix}.
\end{equation}
\end{prp}

The method in~\cite{DDCQMI:kaminaga2025datainformativity} treats~\eqref{apx:LMI} as an LMI in $(\Phi,\alpha)$ and optimizes a control objective. We use the same condition to minimize the volume of the outer approximation:
\begin{equation}\label{prb:relaxed}
\begin{aligned}
    \min_{\substack{C\in\mathbb{R}^{q\times r},\ Q\in\mathbb{S}_{++}^q,\\
    R\in\mathbb{S}_{++}^r,\ \alpha\in\mathbb{R}_+^K}}
        \Psi(Q,R)\quad\st\quad\eqref{apx:LMI},\quad\Phi=\Pi(C,Q,R).
\end{aligned}
\end{equation}
Problem~\eqref{prb:relaxed} replaces containment with a sufficient condition and is therefore conservative. An optimal solution of~\eqref{prb:relaxed} need not solve the original problem~\eqref{prb:prb}. Our goal is to reformulate~\eqref{prb:relaxed} as a convex optimization problem without the computationally demanding LMI constraints.

We address two technical difficulties. 
The first difficulty is the computational cost of imposing~\eqref{apx:LMI}.
We address it by deriving an equivalent condition and characterizing the optimal solutions of~\eqref{prb:relaxed}, thereby reducing the problem to optimization over $\alpha$ alone.
The second difficulty is the nonconvex volume objective.
Schur complements allow us to express the constraints of~\eqref{prb:relaxed} as an LMI in $(C,Q,R,\alpha)$ or in $(L,P,R,\alpha)$, where $L=Q^{-1}C$ and $P=Q^{-1}$. However, both $\Psi(Q,R)=r\logdet Q+q\logdet R$ and $\Psi(P^{-1},R)=-r\logdet P+q\logdet R$ are nonconvex. The next two sections address these difficulties in turn.

\section{PARAMETERIZATION OF OUTER MATRIX ELLIPSOIDS}\label{sec:parametrize}

We now reduce~\eqref{prb:relaxed} to an optimization problem in $\alpha$ alone. We begin with an equivalent form of~\eqref{apx:LMI}.
\begin{lem}\label{lem:LMI_henkan}
Suppose that Assumption~\ref{ass:volume} holds. Set $\bar C\coloneqq\sum_{k=1}^{K}F_kC_kG_k$, and, for $\alpha\in\mathbb{R}_{++}^K$, define
\begin{equation}\label{apx:parametrize}
    Q(\alpha)\coloneqq\sum_{k=1}^{K}\frac{1}{\alpha_k}F_kQ_kF_k^\top,\ 
    R(\alpha)\coloneqq\sum_{k=1}^{K}\alpha_kG_k^\top R_kG_k.
\end{equation}
Then $C\in\mathbb{R}^{q\times r}$, $Q\in\mathbb{S}_{++}^q$, $R\in\mathbb{S}_{++}^r$, and $\alpha\in\mathbb{R}_+^K$ satisfy~\eqref{apx:LMI} if and only if $\alpha\in\mathbb{R}_{++}^K$ and
\begin{equation}\label{prm:LMI_CQR}
    \sqmat{R-R(\alpha)}{\bar C^\top-C^\top}{\bar C-C}{Q-Q(\alpha)}\succeq0.
\end{equation}
\end{lem}
\begin{proof}
We first show that~\eqref{apx:LMI} implies $\alpha\in\mathbb{R}_{++}^K$. The lower-right block of~\eqref{apx:LMI} gives
$-F^\top Q^{-1}F+\diag_{k=1}^{K}(\alpha_kQ_k^{-1})\succeq0$.
Each diagonal block gives $\alpha_kQ_k^{-1}-F_k^\top Q^{-1}F_k\succeq0$. Since $F_k\ne0$ by Assumption~\ref{ass:volume}, the matrix $-F_k^\top Q^{-1}F_k$ has a negative eigenvalue. Thus, $\alpha_k=0$ is impossible. We conclude that $\alpha_k>0$ for every $k$, so $\alpha\in\mathbb{R}_{++}^K$.

Next, we establish the equivalence of~\eqref{apx:LMI} and~\eqref{prm:LMI_CQR} for $\alpha\in\mathbb{R}_{++}^K$. Substituting $\Phi=\Pi(C,Q,R)$ and $\Phi_k=\Pi(C_k,Q_k,R_k)$, $k=1,\dots,K$, into~\eqref{apx:LMI} gives
\begin{equation}\label{prf:LMI1}
    \left[\begin{smallmatrix}
        R-R(\alpha)&0\\0&0
    \end{smallmatrix}\right]
    -\left[\begin{smallmatrix}
        -C^\top\\F^\top
    \end{smallmatrix}\right]Q^{-1}\left[\begin{smallmatrix}
        -C^\top\\F^\top
    \end{smallmatrix}\right]^\top
    +\left[\begin{smallmatrix}
        \Gamma^\top D\Gamma&-\Gamma^\top D\\
        -D\Gamma &D
    \end{smallmatrix}\right]\succeq0,
\end{equation}
where $\Gamma=[G_1^\top C_1^\top\cdots G_K^\top C_K^\top]^\top$ and $D=\diag_{k=1}^{K}(\alpha_kQ_k^{-1})$.
Since $Q\succ0$, the Schur complement yields the following equivalent form of~\eqref{prf:LMI1}:
\begin{equation}\label{prf:LMI2}
    \left[\begin{smallmatrix}
        R-R(\alpha)+\Gamma^\top D\Gamma&-\Gamma^\top D&-C^\top\\
         -D\Gamma&D&F^\top\\
         -C&F&Q
    \end{smallmatrix}\right]
    \succeq0.
\end{equation}
Since $\alpha_k>0$ and $Q_k\succ0$, we have $D\succ0$. Taking the Schur complement of the block $D$ and using $\bar{C}=F\Gamma$ and $Q(\alpha)=FD^{-1}F^\top$ yield~\eqref{prm:LMI_CQR}. All these transformations are reversible, which proves the equivalence.
\end{proof}

The next theorem uses this lemma to characterize the optimal solutions of~\eqref{prb:relaxed}: the center $C$ is fixed, and the shape matrices $(Q,R)$ depend only on $\alpha$.
\begin{thm}\label{thm:parametrize}
Suppose that Assumption~\ref{ass:volume} holds, and let $\bar C$, $Q(\alpha)$, and $R(\alpha)$ be as in Lemma~\ref{lem:LMI_henkan}. If $(C^\star,Q^\star,R^\star,\alpha^\star)$ is an optimal solution of~\eqref{prb:relaxed}, then
\begin{equation}
    C^\star=\bar C,\quad Q^\star=Q(\alpha^\star),\quad R^\star=R(\alpha^\star).
\end{equation}
\end{thm}
\begin{proof}
Take any feasible solution $(C,Q,R,\alpha)$ of~\eqref{prb:relaxed}. By Lemma~\ref{lem:LMI_henkan}, we have $\alpha\in\mathbb{R}_{++}^K$ and~\eqref{prm:LMI_CQR}. The diagonal blocks in~\eqref{prm:LMI_CQR} give $Q\succeq Q(\alpha)$ and $R\succeq R(\alpha)$. Moreover,
\begin{align}
    Q(\alpha)&=F\diag_{k=1}^{K}(\alpha_k^{-1}Q_k)F^\top,\\
    R(\alpha)&=G^\top\diag_{k=1}^{K}(\alpha_kR_k)G.
\end{align}
Since $F$ has full row rank and $G$ has full column rank by Lemma~\ref{lem:FG}, both $Q(\alpha)$ and $R(\alpha)$ are positive definite. Hence, 
$R\otimes Q\succeq R(\alpha)\otimes Q(\alpha)\succ0$.
Since the log-determinant is monotone with respect to the L\"owner order, it follows that $\Psi(Q,R)\ge\Psi(Q(\alpha),R(\alpha))$, with equality if and only if $Q=Q(\alpha)$ and $R=R(\alpha)$.

Next, take any $\alpha\in\mathbb{R}_{++}^K$ and set $(C,Q,R)=(\bar C,Q(\alpha),R(\alpha))$. As above, Lemma~\ref{lem:FG} ensures that $Q$ and $R$ are positive definite.
Then,~\eqref{prm:LMI_CQR} is satisfied, meaning that $(\bar C,Q(\alpha),R(\alpha),\alpha)$ is feasible for~\eqref{prb:relaxed}.

Apply the preceding argument to an optimal solution
$(C^\star,Q^\star,R^\star,\alpha^\star)$.
The feasible point
$(\bar C,Q(\alpha^\star),R(\alpha^\star),\alpha^\star)$
has no larger objective value, so optimality requires equality.
Hence, $Q^\star=Q(\alpha^\star)$ and $R^\star=R(\alpha^\star)$.
Substituting these equalities into~\eqref{prm:LMI_CQR} gives $C^\star=\bar C$.
\end{proof}

By Theorem~\ref{thm:parametrize}, problem~\eqref{prb:relaxed} reduces to
\begin{equation}\label{prb:alpha}
    \min_{\alpha\in\mathbb{R}_{++}^K}
    f(\alpha)\coloneqq\Psi(Q(\alpha),R(\alpha)).
\end{equation}
If $\alpha^\star$ solves~\eqref{prb:alpha}, then $(\bar C,Q(\alpha^\star),R(\alpha^\star),\alpha^\star)$ solves~\eqref{prb:relaxed}. Note that the objective is invariant under scaling. For all $a>0$, it holds that
\begin{equation}\label{obj_scale_invariant}
    f(a\alpha)=f(\alpha).
\end{equation}
This reflects the scaling freedom~\eqref{e:scale_invariant} of matrix ellipsoids. Since $Q(a\alpha)=a^{-1}Q(\alpha)$ and $R(a\alpha)=aR(\alpha)$, the parameters $\alpha$ and $a\alpha$ describe the same outer matrix ellipsoid.

\emph{Remark (Relation to the vector-valued case):} For ordinary ellipsoids,~\eqref{apx:parametrize} reduces to the parameterization in~\cite{ellip:durieu_multi-input_2001}. To see this, let $r=1$ and use the scaling freedom~\eqref{e:scale_invariant} to normalize $G_k^\top R_kG_k=1$ and $R(\alpha)=1$. Then~\eqref{apx:parametrize} becomes
\begin{equation}\label{apx:parametrize_vec}
    Q(\alpha)=\sum_{k=1}^{K}\frac{1}{\alpha_k}F_kQ_kF_k^\top,\ 
    \sum_{k=1}^{K}\alpha_k=1,
\end{equation}
which is equivalent to the representation in~\cite[Theorem 4.2]{ellip:durieu_multi-input_2001}. That work derives the representation using support functions. Our derivation instead follows from the LMI condition in Proposition~\ref{prp:LMI}, which extends the condition in~\cite{boyd1994linear} to matrix variables. In the vector-valued case, our result therefore establishes equivalence between volume minimization based on the S-procedure and optimization over the parameterized family. With~\eqref{apx:parametrize_vec}, problem~\eqref{prb:alpha} becomes
\begin{equation}\label{prb:vec}
    \min_{\alpha\in\mathbb{R}_{++}^K}
        \log\det\!\left(\sum_{k=1}^{K}
              \frac{F_kQ_kF_k^\top}{\alpha_k}\right)\quad\st\quad\sum_{k=1}^{K}\alpha_k=1,
\end{equation}
whose convexity was established in~\cite[Theorem 4.3]{ellip:durieu_multi-input_2001}.

\section{OPTIMAL MATRIX ELLIPSOIDS}\label{sec:opt}
\subsection{Equivalent Convex Reformulation}

Problem~\eqref{prb:alpha} has only the scalar variables $\alpha\in\mathbb{R}_{++}^K$ and no LMI constraints. Its objective, however, is generally nonconvex. For example, let $q=r=2$ and $D=\diag(10,16)$, and consider
\begin{equation}
    Q(\alpha)=\frac{I_2}{\alpha_1}+\frac{D}{\alpha_2},\quad 
    R(\alpha)=\alpha_1D+\alpha_2I_2.
\end{equation}
Then
$f(\alpha)=4\log\frac{(10\alpha_1+\alpha_2)(16\alpha_1+\alpha_2)}{\alpha_1\alpha_2}$,
and $\frac{d^2}{dt^2}f(t,1-t)\rvert_{t=1/3}=-\frac14<0$.
Thus, $f$ is not convex.

A logarithmic change of variables resolves this nonconvexity. For $\beta\in\mathbb{R}^K$, let $e^\beta\in\mathbb{R}^K$ denote the vector whose $k$th entry is $e^{\beta_k}$. Substituting $\alpha=e^\beta$ transforms~\eqref{prb:alpha} into an unconstrained minimization problem with objective $g(\beta)\coloneqq f(e^\beta)$. This objective can be written as
\begin{align}
    g(\beta)=q\logdet\!\left(\sum_{k=1}^{K}e^{\beta_k}\widetilde R_k\right)\!+r\logdet\!\left(\sum_{k=1}^{K}e^{-\beta_k}\widetilde Q_k\right),
\end{align}
where $\widetilde Q_k=F_kQ_kF_k^\top$ and $\widetilde R_k=G_k^\top R_kG_k$. If $q=r=1$, then $g$ is a sum of log-sum-exp functions and is therefore convex~\cite{boyd2004convex}. We establish convexity for general $q$ and $r$.

\begin{thm}\label{thm:logdet_convex}
Suppose that Assumption~\ref{ass:volume} holds. Then $g(\beta)$ is convex, and its gradient is given by
\begin{equation}
    \frac{\partial g}{\partial\beta_k}
    ={}q\tr(R^{-1}\widetilde R_k)e^{\beta_k}-r\tr(Q^{-1}\widetilde Q_k)e^{-\beta_k},
\end{equation}
where $Q$ and $R$ abbreviate $Q(e^\beta)$ and $R(e^\beta)$. Its Hessian $H$ is given by $H=\Xi-\Omega$ with
\begin{align}
    \Xi={}&\diag_{k=1}^{K}\bigl(q\tr(R^{-1}\widetilde R_k)e^{\beta_k}+r\tr(Q^{-1}\widetilde Q_k)e^{-\beta_k}\bigr),\\
    \Omega_{k,l}={}&q\tr(R^{-1}\widetilde R_kR^{-1}\widetilde R_l)
                        e^{\beta_k}e^{\beta_l}\\
        &+r\tr(Q^{-1}\widetilde Q_kQ^{-1}\widetilde Q_l)
                        e^{-\beta_k}e^{-\beta_l}.
\end{align}
\end{thm}
\begin{proof}
To analyze both terms of $g$ at once, consider
$g^M(\beta)\coloneqq\logdet M$, where
$M=\sum_{k=1}^{K}M_ke^{s\beta_k}$ and $s\in\{-1,1\}$.
Assume that $M_k\succeq0$ for every $k$ and
$\sum_{k=1}^{K}M_k\succ0$, so that $M\succ0$ for every $\beta$.
Standard differentiation formulas~\cite[Sec. 6.3]{ben2001lectures} give
\begin{align}
    \frac{\partial g^M}{\partial\beta_k}
        =&{}\tr(M^{-1}M_k)se^{s\beta_k},\\
    H^M_{k,l}
        =&-\tr(M^{-1}M_lM^{-1}M_k)e^{s\beta_k}e^{s\beta_l}\\
        &+\delta_{k,l}\tr(M^{-1}M_k)e^{s\beta_k},
\end{align}
where $H^M$ is the Hessian of $g^M$, and $\delta_{k,l}$ equals $1$ if $k=l$ and $0$ otherwise. In deriving the Hessian, we used $s^2=1$.

To prove convexity of $g^M$, we show that $v^\top H^Mv\ge0$ for every $v\in\mathbb{R}^K$. Define $B_{k,l}=\tr(M^{-1}M_lM^{-1}M_k)e^{s\beta_k}e^{s\beta_l}$. Since 
$\tr(M^{-1}M_k)e^{s\beta_k}=\tr(M^{-1}M_kM^{-1}M)e^{s\beta_k}=\sum_{l=1}^{K}B_{k,l}$,
the quadratic form can be written as
{\small\begin{equation}\label{prf:hessian}
\begin{split}
    2v^\top\! H^Mv
        ={}&\textstyle{2\sum_{k=1}^{K}\!v_k^2\tr(M^{-1}M_k)e^{s\beta_k}\!-\!2\sum_{k=1}^{K}\!\sum_{l=1}^{K}\!v_kv_lB_{k,l}}\\
        ={}&\textstyle{\sum_{k=1}^{K}\sum_{l=1}^{K}(v_k-v_l)^2B_{k,l}.}
\end{split}
\end{equation}}Because $M\succ0$ and $M_k\succeq0$, every $B_{k,l}$ is nonnegative. Hence, $v^\top H^Mv\ge0$.

As shown in the proof of Theorem~\ref{thm:parametrize}, Assumption~\ref{ass:volume} ensures that $Q(e^\beta)$ and $R(e^\beta)$ are positive definite for every $\beta\in\mathbb{R}^K$. Applying the results for $g^M$ to each term of $g$ gives the claimed convexity, gradient, and Hessian.
\end{proof}

We can now minimize $g$ using standard convex optimization methods. The scaling freedom~\eqref{obj_scale_invariant} becomes $g(\beta+(\log a)\mathbf{1})=g(\beta)$ for every $a>0$. Thus, we can impose the normalization $\mathbf{1}^\top\beta=0$ without changing the objective value and solve
\begin{equation}\label{prb:beta}
    \min_{\beta\in\mathbb{R}^K}g(\beta)
    \quad\st\quad\mathbf{1}^\top\beta=0.
\end{equation}
The next proposition shows that $g$ is strictly convex on the normalized space and that its gradient lies in the same space. These properties are useful for numerical optimization.

\begin{prp}\label{prp:strict_cvx}
Suppose that Assumption~\ref{ass:volume} holds, and define $\mathcal{F}_\beta\coloneqq\{\beta\in\mathbb{R}^K\mid\mathbf{1}^\top\beta=0\}$. Then $\nabla g(\beta)\in\mathcal{F}_\beta$, and $g$ is strictly convex on $\mathcal{F}_\beta$.
\end{prp}
\begin{proof}
By~\eqref{obj_scale_invariant}, $g(\beta+b\mathbf{1})=g(\beta)$ for every $b\in\mathbb{R}$. Differentiating with respect to $b$ gives $\mathbf{1}^\top\nabla g(\beta)=0$, so $\nabla g(\beta)\in\mathcal{F}_\beta$.

We next prove strict convexity by examining the Hessian $H$ in Theorem~\ref{thm:logdet_convex}. A calculation analogous to~\eqref{prf:hessian} yields
\begin{equation}\label{prf:quadform}
    \textstyle{2v^\top Hv=\sum_{k=1}^{K}\sum_{l=1}^{K}(v_k-v_l)^2\Omega_{k,l}.}
\end{equation}
The trace properties also give
\begin{equation}
\begin{split}
    \Omega_{k,l}={}&q\left\lVert
        R_k^{1/2}G_kR^{-1}G_l^\top R_l^{1/2}
        \right\rVert_F^2 e^{\beta_k}e^{\beta_l}\\
        &+r\left\lVert
        Q_k^{1/2}F_k^\top Q^{-1}F_lQ_l^{1/2}
        \right\rVert_F^2 e^{-\beta_k}e^{-\beta_l},
\end{split}
\end{equation}
where $\norm{M}_F$ denotes the Frobenius norm of $M$. Since $Q_k$ and $R_k$ are positive definite for every $k$,
\begin{equation}\label{prf:psi_gf}
    \Omega_{k,l}=0\ \Rightarrow\
    G_kR^{-1}G_l^\top=0\text{ and }
    F_k^\top Q^{-1}F_l=0.
\end{equation}

Take a nonzero $v\in\mathcal{F}_\beta$ and suppose that $v^\top Hv=0$. By~\eqref{prf:quadform}, either $v_k=v_l$ or $\Omega_{k,l}=0$ for every pair $(k,l)$. Since $v\in\mathcal{F}_\beta$ and $v\ne0$, there exist indices $i$ and $j$ such that $v_i>0>v_j$. For each $k$, we then have $v_k\ne v_i$ or $v_k\ne v_j$, and hence $\Omega_{k,i}=0$ or $\Omega_{k,j}=0$.
It follows from~\eqref{prf:psi_gf} that, for every $k$, $F_k^\top Q^{-1}F_i=0$ or $G_kR^{-1}G_j^\top=0$.
For $y\in\im(Q^{-1}F_i)$ and $z\in\im(R^{-1}G_j^\top)$,
either $F_k^\top y=0$ or $G_kz=0$ holds for every $k$.
Hence, every $X=\sum_{k=1}^{K}F_kX_kG_k\in\mathcal{T}$,
with $X_k\in\el_k$, satisfies
\begin{equation}
    \textstyle{y^\top Xz=\sum_{k=1}^{K}(F_k^\top y)^\top X_k(G_kz)=0.}
\end{equation}
Since $F_i$ and $G_j$ are nonzero by Assumption~\ref{ass:volume}, we can choose $y\ne0$ and $z\ne0$. Then $\mathcal{T}$ lies in the proper linear subspace $\{X\in\mathbb{R}^{q\times r}\mid y^\top Xz=0\}$, contradicting the nonempty interior assumption. Therefore, $v^\top Hv>0$, which proves strict convexity.
\end{proof}

For gradient descent, consider the update $\betat{t+1}=\betat{t}-\tau\nabla g(\betat{t})$ with $\betat{0}\in\mathcal{F}_\beta$ and $\tau>0$. Proposition~\ref{prp:strict_cvx} gives $\nabla g(\betat{t})\in\mathcal{F}_\beta$, so every iterate remains feasible. For Newton's method, we compute a correction $d\in\mathcal{F}_\beta$ by solving $Hd=\nabla g(\beta)$. Since $H$ is positive definite on $\mathcal{F}_\beta$, this equation has a unique solution in that subspace.

\subsection{Comparison of Computational Complexity}

We compare the cost per iteration of two methods: an interior-point method that directly imposes~\eqref{apx:LMI}, and Newton's method applied to~\eqref{prb:beta}. For this comparison, we view the interior-point method as Newton's method applied to a barrier formulation of the constrained problem.

For~\eqref{prb:beta}, each Newton iteration requires solving $Hd=\nabla g(\beta)$, where $H$ is the Hessian in Theorem~\ref{thm:logdet_convex}. This typically takes $O(K^3)$ operations. We must also evaluate $H=\Xi-\Omega$ and $\nabla g$; forming $\Omega$ dominates this cost. We first compute $R^{-1}\widetilde R_k$ and $Q^{-1}\widetilde Q_k$ for every $k$ in $O(K(q^3+r^3))$ operations. Each $\Omega_{k,l}$ is then obtained from traces of products of these matrices in $O(q^2+r^2)$ operations. Assembling all entries of $\Omega$ therefore takes $O(K^2(q^2+r^2))$ operations. The total cost per Newton iteration is
\begin{equation}\label{eq:newton_complexity}
    O\!\left(K(q^3+r^3)+K^2(q^2+r^2)+K^3\right).
\end{equation}

Next, consider an interior-point method that treats~\eqref{apx:LMI} as an LMI in $(\Phi,\alpha)$. The method in~\cite[Sec. VI]{DDCQMI:kaminaga2025datainformativity} uses this formulation to minimize control objectives subject to~\eqref{apx:LMI} and additional constraints. For simplicity, we omit the cost of evaluating the objective and handling those additional constraints. There are $N_{\rm var}=(q+r)(q+r+1)/2+K$ scalar decision variables in~\eqref{apx:LMI}, and the LMI has size $p\times p$, where $p=r+\sum_{k=1}^{K}q_k$. An interior-point method for this constraint has a per-iteration cost of
\begin{equation}
    O\!\left(N_{\rm var}p^3+N_{\rm var}^2p^2+N_{\rm var}^3\right)
\end{equation}\cite{boyd2004convex,ben2001lectures}. 
This cost is higher than~\eqref{eq:newton_complexity}. 
With $q$, $r$, and $\bar q=\frac{1}{K}\sum_kq_k$ fixed,
the cost per iteration is $O(K^4)$ for the LMI formulation
and $O(K^3)$ for the proposed formulation.
With $K$ fixed and $p=O(q+r)$, the corresponding costs are
$O(q^6+r^6)$ and $O(q^3+r^3)$, respectively.

\section{NUMERICAL COMPARISONS}\label{sec:numerical}

We evaluate approximation quality and computation time by comparing two methods. The proposed method applies Newton's method to~\eqref{prb:beta}; the baseline uses the LMI in Proposition~\ref{prp:LMI}. We solved the semidefinite programs (SDPs) using the MATLAB function \texttt{mincx}~\cite{cntr:Gahinet_LMIcontrolTool}. All experiments were run on a PC with an Intel Xeon W-2295 processor (3.00~GHz) and 64~GB of RAM. We executed the trials and methods sequentially.

The reported computation times include all processing required by each method. For the LMI-based method, this includes constructing the constraints and numerically checking the returned solution. In the figures, points indicate medians, and error bars indicate the first and third quartiles. The timing statistics include all trials, even those in which the solution was not accepted.

\subsection{Comparison for Volume Minimization}

We first compare the methods on randomly generated outer approximation problems. The baseline fixes $R=I_r$ and $C=\bar C$ in~\eqref{prb:relaxed}, introduces $P=Q^{-1}$, and solves
\begin{equation}\label{prb:for_compare}
    \min_{P\in\mathbb{S}_{++}^q,\ \alpha\in\mathbb{R}_+^K}
        -\logdet P\quad\st\quad\eqref{apx:LMI},\ 
                          \Phi=\Pi(\bar C,P^{-1},I_r).
\end{equation}
With $C$ and $R$ fixed, this problem is convex in $(P,\alpha)$.\footnote{Because \texttt{mincx} accepts only linear objectives, our implementation uses an equivalent reformulation that recursively introduces $2\times2$ LMIs.}

We varied the dimensions $(q,r)$ and the number of matrix ellipsoids $K$, setting $F_k=I_q$ and $G_k=I_r$ for every $k$. For each problem, we set $C_k=0$ and generated $Q_k=M_{k}^\top M_{k}$ and $R_k=N_k^\top N_k$ by drawing the entries of $M_k\in\mathbb{R}^{q\times q}$ and $N_k\in\mathbb{R}^{r\times r}$ from the standard normal distribution. We regenerated any matrices that were not positive definite.

\begin{figure}[t]
    \centering
    \includegraphics[width=\columnwidth]{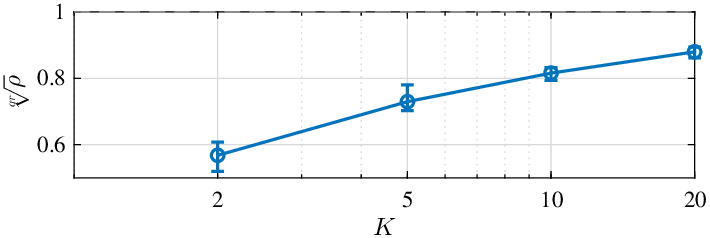}
    \caption{Normalized volume ratio $\sqrt[qr]{\rho}$ as $K$ varies, with $q=r=5$.}
    \label{fig:vol_K}
\end{figure}

Both methods successfully solved all problem instances in this subsection.
First, we fixed $q=r=5$ and ran $50$ trials for each $K\in\{2,5,10,20\}$. Both methods used the same problem in each trial. 
Let $\alpha^\star$ be the weights corresponding to a solution of~\eqref{prb:beta}, and let $P^\star$ solve~\eqref{prb:for_compare}. We define the volume ratio of $\el(\bar C,Q(\alpha^\star),R(\alpha^\star))$ to $\el(\bar C,P^{\star-1},I_r)$ as
\begin{equation}
    \rho=\sqrt{\det\!\left(R(\alpha^\star)\otimes Q(\alpha^\star)\right)}
    \,/
    \sqrt{\det\!\left(I_r\otimes{P^{\star-1}}\right)}.
\end{equation}
To account for the ambient dimension $qr$, we report the normalized
volume ratio $\sqrt[qr]{\rho}$.
Fig.~\ref{fig:vol_K} shows that its median and quartiles increase with $K$, but both the median and the third quartile remain below $1$ for every tested value of $K$. Thus, in this experiment, the proposed method tends to produce outer approximations with smaller volumes than the baseline. 
This difference reflects the shape restriction $R=I_r$ in the baseline.
The proposed method solves the full problem~\eqref{prb:relaxed}
through its equivalent formulation~\eqref{prb:beta}.

\begin{figure}[t]
    \centering
    \includegraphics[width=\columnwidth]{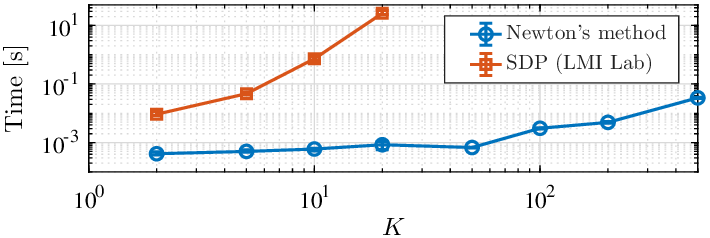}
    \caption{Computation time as $K$ varies, with $q=r=5$. Newton's method solves~\eqref{prb:beta}, whereas SDP (LMI Lab) solves~\eqref{prb:for_compare}.}
    \label{fig:Time_K}
\end{figure}

Fig.~\ref{fig:Time_K} shows the computation times. We also timed the proposed method in $50$ trials for each $K\in\{50,100,200,500\}$. Its median computation time is lower at every $K$ for which both methods were tested, and the difference grows with $K$. Even at $K=500$, the median for the proposed method remains below $0.1$~s.

\begin{figure}[t]
    \centering
    \includegraphics[width=\columnwidth]{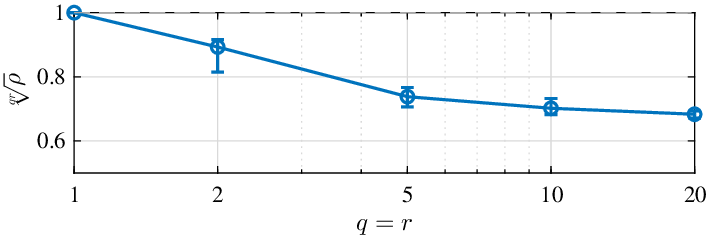}
    \caption{Normalized volume ratio $\sqrt[qr]{\rho}$ as the dimension varies, with $K=5$ and $q=r$.}
    \label{fig:vol_qr}
\end{figure}

\begin{figure}[t]
    \centering
    \includegraphics[width=\columnwidth]{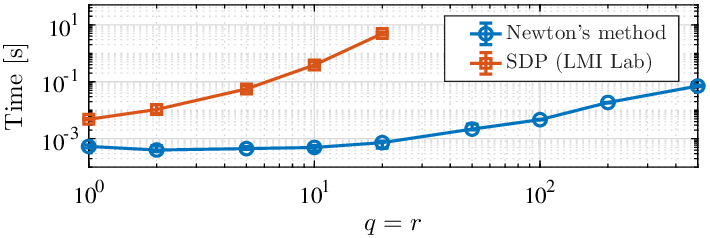}
    \caption{Computation time as the dimension varies, with $K=5$ and $q=r$.}
    \label{fig:Time_qr}
\end{figure}

Next, we fixed $K=5$ and $q=r$ and ran $50$ trials for each $q\in\{1,2,5,10,20\}$. We also timed the proposed method in $50$ trials for each $q\in\{50,100,200,500\}$. Figs.~\ref{fig:vol_qr} and~\ref{fig:Time_qr} show the normalized volume ratio $\sqrt[qr]{\rho}$ and computation time, respectively. In Fig.~\ref{fig:vol_qr}, the ratio is below $1$ for all tested settings with $q=r\ge2$. Its median decreases as the dimension increases, showing a greater improvement in approximation volume. For $q=r=1$, the ratio equals $1$ up to numerical error. This agrees with~\eqref{prb:vec}: fixing $R=1$ entails no loss of generality when $r=1$. Fig.~\ref{fig:Time_qr} also shows lower median computation times for the proposed method than for the baseline. The median for the proposed method remains below $0.1$~s even at $q=r=500$.

\subsection{Application to Data-Driven Control}

We next apply the proposed method to data-driven control under the noise models in Section~\ref{subsec:cnt}. For these models, the method in~\cite{DDCQMI:kaminaga2025datainformativity} jointly determines an outer approximation and a controller by solving
\begin{equation}\label{exp:codesign}
\begin{aligned}
    \min_{\theta\in\Theta,\ \Phi\in\mathbb{S}^{q+r},\ \alpha\in\mathbb{R}_+^K}
    c^\top \theta\quad\st\quad\eqref{apx:LMI},\  h(\theta,\Phi)\succeq0,
\end{aligned}
\end{equation}
where $\theta$ is the vector of controller design parameters, $\Theta$ is its domain, and $c$ is a vector of the same dimension as $\theta$. The matrix-valued function $h(\theta,\Phi)$ is affine in $(\theta,\Phi)$.

Our approach separates outer approximation from controller design. First, we solve~\eqref{prb:beta} for $\beta^\star$, recover the corresponding weights $\alpha^\star$, and set $\Phi_{\rm fixed}\coloneqq\Pi(\bar C,Q(\alpha^\star),R(\alpha^\star))$. This matrix represents the outer approximation. With its shape fixed, we solve
\begin{equation}\label{exp:twostep}
    \min_{\theta\in\Theta,\ \eta\ge0} c^\top \theta\quad\st \quad h(\theta,\eta\Phi_{\rm fixed})\succeq0,
\end{equation}
where the scalar $\eta$ allows the scaling of the QMI coefficient matrix $\Phi$ to be optimized during controller design.
We compare this two-stage approach, which solves~\eqref{prb:beta} followed by~\eqref{exp:twostep}, with the joint formulation~\eqref{exp:codesign} in terms of guaranteed control performance and computation time.

For the unknown true system in Section~\ref{subsec:cnt}, we used a model derived from the inverted-pendulum model in~\cite[(14)]{Con:chalupa2008}. We discretized it with a sampling period of $0.04$~s and retained three states: velocity, pendulum angle, and angular velocity. The resulting system matrices are
\begin{align}
    A^*=\left[\begin{smallmatrix*}[r]
        0.9384&-0.0253&-0.0004\\
        0.0022& 1.0150& 0.0401\\
        0.1108& 0.7513& 1.0117
    \end{smallmatrix*}\right],\ 
    B^*=\left[\begin{smallmatrix*}[r]
        0.0095\\-0.0003\\-0.0171
    \end{smallmatrix*}\right].
\end{align}
The state and input dimensions are $n=3$ and $m=1$. We added measurement noise to all three data matrices $(X_+,X,U)$ and assumed a known common amplitude bound: $\Delta_{k,l}^2\le\epsilon^2$ for every $(k,l)$, with $\epsilon=2\times10^{-4}$. This is the element-wise bound in~\ref{cnt:elementwise_bound} of Section~\ref{subsec:cnt}, and thus can be described by a Minkowski sum. For controller design in~\eqref{exp:twostep}, we used the $H_2$ performance condition in~\cite[(35)]{DDCQMI:kaminaga2025datainformativity}. 
We evaluate the closed-loop $H_2$ performance from the disturbance
$w\in\mathbb{R}^n$ to the output
$z=[x^\top\ 10^{-5}u^\top]^\top$,
with dynamics $x_+=A^*x+B^*u+w$.

In each trial, we sampled the entries of the noise-free input and state data $U^*$ and $X^*$ independently from the standard normal distribution and computed $X_+^*=A^*X^*+B^*U^*$ as in~\eqref{cnt:cleandata}. We then sampled the noise entries independently from the uniform distribution on $[-\epsilon,\epsilon]$ and formed the measured data using~\eqref{cnt:measured_data}. We ran $50$ trials for each data length $T\in\{5,7,9,11,13\}$. Both methods used the same measured data and noise bound in each trial.

We report computation time and the performance ratio $\gamma/\gamma^\star$, where $\gamma$ is the guaranteed performance and $\gamma^\star$ is the optimal $H_2$ performance achievable with the true model known. 
Performance statistics are computed separately for each method using
the trials in which it returned a feasible solution.
Timing statistics include all trials.
The reported time for the proposed method includes both the approximation step~\eqref{prb:beta} and the controller design step~\eqref{exp:twostep}.

\begin{figure}[t]
    \centering
    \includegraphics[width=\columnwidth]{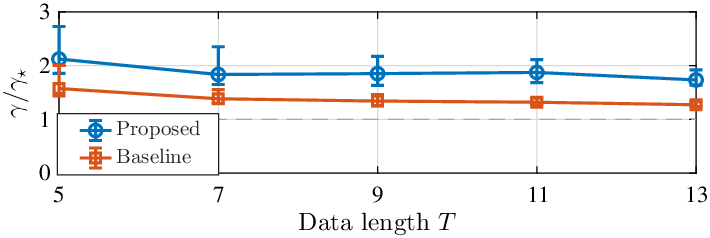}
    \caption{Ratio to the optimal performance. Proposed first computes the approximation using the proposed method and then solves~\eqref{exp:twostep}; Baseline solves~\eqref{exp:codesign}.}
    \label{fig:h2_ratio}
\end{figure}

\begin{figure}[t]
    \centering
    \includegraphics[width=\columnwidth]{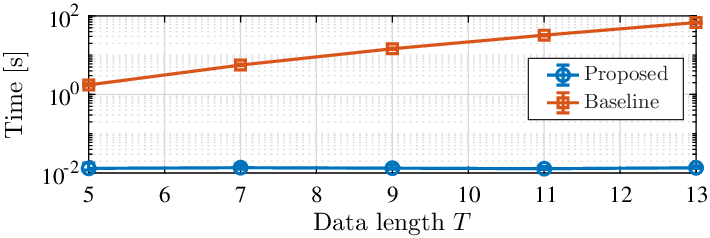}
    \caption{Computation time as the data length $T$ varies.}
    \label{fig:cnt_time}
\end{figure}

Figs.~\ref{fig:h2_ratio} and~\ref{fig:cnt_time} show the results, with the two-stage method labeled Proposed and the joint formulation labeled Baseline. In Fig.~\ref{fig:h2_ratio}, Baseline has a lower performance ratio than Proposed for every $T$. The proposed method first minimizes the volume of the outer approximation to the noise set. This approximation need not minimize the guaranteed upper bound on the $H_2$ performance. Baseline determines the approximation and controller jointly, so it can adjust the approximation shape to suit the control objective. This additional design freedom explains the difference.

In Fig.~\ref{fig:cnt_time}, the median computation time of the proposed method stays on the order of $10^{-2}$~s across the tested values of $T$. The computation time of Baseline increases substantially with $T$, and its median is at least two orders of magnitude higher at every $T$. Fixing the approximation shape in advance reduces the number of decision variables and constraints in the controller design step. In this setting, the resulting reduction in computation time comes at the cost of more conservative performance guarantees.

\section{CONCLUSIONS}\label{sec:conclusion}

We have derived an equivalent convex formulation for minimizing the volume of an outer matrix ellipsoid of a Minkowski sum under LMI conditions based on the matrix S-procedure.
By expressing the optimal center and shape in terms of scalar parameters, we eliminate the matrix variables and LMI constraints.
A logarithmic change of variables then makes the objective convex.
We also prove strict convexity on the normalized space and derive the gradient and Hessian, which allow us to apply Newton's method and analyze its computational cost. 
Numerical experiments show that the proposed method is substantially faster than the conventional LMI-based method.
In the control example, computing the approximation in advance substantially reduced
computation time but produced more conservative performance guarantees.
It therefore offers an option when computational savings justify
this additional conservatism.

\bibliographystyle{IEEEtran}
\bibliography{IEEEabrv,Control,Robust_DDC,DDC,ellipsoid,mathematics}

\end{document}